\documentclass[12pt,twoside]{article}
\usepackage{amssymb}
\usepackage[mathscr]{eucal}
\usepackage{eufrak}
\usepackage{amsmath}
\usepackage{mathrsfs}
\usepackage[colorlinks=true,
   urlcolor=blue,           filecolor=green,      
   citecolor=green,      
   linkcolor=red,           bookmarks=true,
  unicode,
   plainpages=false,   ]{hyperref}
\usepackage{color}

\def\intl#1{\int\limits_{#1}}
\def\intll#1#2{\int\limits_{#1}^{#2}}
\def\dm{|\hskip-0.05cm|}
\def\OO{\Omega}
\def\displ{\displaystyle}
\def\VSE{\vspace{6pt}\\&\displ }
\def\VS{\vspace{6pt}\\\displ }
\def\rf#1{{\rm(\ref{#1})}}

\def\chiu{\hfill$\displaystyle\vspace{4pt}
\underset\Box\null$\par}
\def\Pr{{\bf Proof. }}
\def\O{\Omega}

\def\R{\Bbb R}
\def\N{\Bbb N}

\def\à{\`{a}}
\def\è{\`{e}}
\def\ì{\`{i}}
\def\ù{\`{u}}
\def\ò{\`{o}}
\def\é{\'{e}}

\def\vep{\varepsilon}

\def\be{\begin{equation}}
\def\ba{\begin{array}}
\def\ea{\end{array}}
\def\ee{\end{equation}}
\def\vs1{\vspace{1ex}}

\def\ov{\overline}

\font\sc=cmcsc10
\title{Some remarks on the partial regularity of a suitable weak solution to the Navier-Stokes Cauchy problem}
\begin{document}
\markboth{\footnotesize\it  F. Crispo and P.
Maremonti} {\footnotesize\it On the partial regularity of a suitable weak solution to the Navier-Stokes Cauchy problem}\author{\sc F. Crispo and P. Maremonti
\thanks{
Dipartimento di Matematica e Fisica, 
Universit\`{a} degli Studi della Campania ``Luigi Vanvitelli", via Vivaldi 43, 81100 Caserta,
 Italy.  Email: 
francesca.crispo@unicampania.it;
paolo.maremonti@unicampania.it}
\\{\footnotesize\it To Professor Vsevolod Alekseevich  Solonnikov on his 85th birthday}\\{\footnotesize\it with the best compliments}}

\date{}

\maketitle \noindent{\bf Abstract} - {\small The aim of the paper is to investigate on some questions of local regularity of a suitable weak solution to the Navier-Stokes Cauchy problem. The  results are obtained in the wake of the ones, well known, by Caffarelli-Kohn-Nirenberg.} \vskip 0.2cm
 \par\noindent{\small Keywords: Navier-Stokes equations, suitable weak solutions, partial regularity. }
  \par\noindent{\small  
  AMS Subject Classifications: 35Q30, 35B65, 76D03.}  
 \par\noindent
 \vskip -0.7true cm\noindent
\newcommand{\red}{\protect\bf}
\renewcommand\refname{\centerline
{\red {\normalsize \bf References}}}
\newtheorem{ass}
{\bf Assumption}[section]
\newtheorem{defi}
{\bf Definition}[section]
\newtheorem{tho}
{\bf Theorem}[section]
\newtheorem{rem}
{\sc Remark}[section]
\newtheorem{lemma}
{\bf Lemma}[section]
\newtheorem{coro}
{\bf Corollary}[section]
\newtheorem{prop}
{\bf Proposition}[section]
\renewcommand{\theequation}{\thesection.\arabic{equation}}
\setcounter{section}{0}
\numberwithin{equation}{section}
\section{Introduction}\label{intro} We deal with the Navier-Stokes Cauchy problem     \be\label{NS}\ba{l}u_t+u\cdot
\nabla u+\nabla\pi_u=\Delta
u,\;\nabla\cdot
u=0,\mbox{ in }(0,T)\times\R^3,\\
u(0,x)=u_0(x)\mbox{ on
}\{0\}\times\R^3.\ea\ee   In system
\rf{NS} $u$ is the kinetic field,
$\pi_u$ is the pressure field,
 $u_t:=
\frac\partial{\partial t}u$  and
 $u\cdot\nabla u:=
u_k\frac\partial{\partial x_k}u$. We investigate on the partial regularity of a suitable weak solution, and we detect a new sufficient condition for the existence of a regular solution.
  Our results are in the wake of the ones obtained in \cite{CKN} and, for small data, in \cite{CMRR}. As in \cite{CMSD,CMRR,MPSC}, our study attempts to highlight what is possible to obtain, without extra condition, in the setting of the $L^2$-theory.  In this connection, although it is not our chief aim, we like to point out that  our results could lead to a sort of structure theorem in the space-time cylinder. To be more precise in the claim  we recall the well known Leray's structure theorem related to a weak solution. Leray's theorem claims that there exist an interval of regularity of the kind $(\theta,\infty)$ and a sequence of intervals of regularity  included in $(0,\theta)$ whose   complementary set on $(0,\theta)$ is a set of zero $\frac12$-Hausdorff measure. {\it Mutatis mutandis},  the results of \cite{CKN} (see below Theorem\,\ref{DCKN-C}) and of this note give a sort of structure theorem for a suitable weak solution related to the Cauchy problem. More precisely,   under a suitable assumption for the initial data, in Theorem\,\ref{DCKN-C} it is proved that a suitable weak solution is regular for all $t>0$ in the exterior of a ball with radius $R_0$. In this note we prove that,   almost everywhere, a point $(t,x)\in (0,\theta)\times B(R_0)$ is the center of a parabolic neighborhood of regularity for a suitable weak solution. Hence in $(0,\theta)\times B(R_0)$ there is at most a sequence of open sets of regularity, whose complementary set in $(0,\theta)\times B(R_0)$ has at most zero $1$-Hausdorff measure.
\par  To
better state the details of our main results, we split the introduction in two short subsections. In the first one we recall
some definitions and notation  following the ones in \cite{CKN}. Then we recall two fundamental regularity results obtained in \cite{CKN}, and, with an alternative proof, in \cite{V}, and their consequences.  In the second subsection we give the statement of our results.
\subsection{\label{CKNR}\bf Suitable weak solutions}
We start by recalling the following:
  \begin{defi}\label{WS}{\sl  Let $u_0\in J^2(\R^3)$. A pair $(u,\pi_u)$, such that $u:(0,\infty)\times\R^3 \to\R^3$ and $\pi_u:(0,\infty)\times\R^3 \to \R$, is said a weak solution to problem {\rm\rf{NS}} if
\begin{itemize}\item [i)] for all $T>0$,
 $u\in  L^2(0,T; J^{1,2}
(\R^3 ))$ and $\pi_u\in L^\frac53((0,T)\times\R^3)$
\item
[ii)] $\displ\lim_{t\to0}\dm
u(t)-u_0\dm_2=0$,\item[iii)]
for all $t,s\in(0,T)$,  the pair $(u,\pi_u)$
satisfies the  equation:
\newline \centerline{$\displ\intll
st\Big[(u,\varphi_\tau)-(\nabla
u,\nabla
\varphi)+(u\cdot\nabla\varphi,u)+(\pi_u,\nabla
\cdot\varphi)\Big]d\tau+(u(s),\varphi
(s))=(u(t),\varphi(t))$,}
\newline\centerline{
 for all $\varphi\in C^1_0([0,T)\times\R^3 )$.}
\end{itemize}}\end{defi} In \cite{CKN} in order to investigate on the regularity of a weak solution it is introduced an     energy relation having a local character:
\begin{defi}\label{SWS}{\sl A pair $(u,\pi_u)$ is said a suitable weak solution if it is a weak solution in the sense of the Definition\,\ref{WS} and, moreover,
\be\label{SEI}\ba{l}\displ\intl{\R^3}|u(t)|^2\phi(t)dx+2\intll \sigma t\intl
{\R^3}|\nabla u|^2\phi\, dxd\tau\leq \intl{\R^3}|u(\sigma)|^2\phi(\sigma)dx\VS\hskip 2,5cm+\intll \sigma t\intl{\R^3}|u|^2(\phi_\tau+\Delta\phi)dxd\tau+\intll \sigma t\intl{\R^3}(|u|^2+2\pi_u)
u\cdot\nabla\phi dxd\tau,\ea\ee for all $t\geq \sigma$, for $\sigma=0$ and a.e. in $\sigma\geq 0$, and for all 
nonnegative $\phi\in C_0^\infty(\R\times\R^3)$.}\end{defi} 
 In \cite{CKN} and \cite{scheffer1} the following existence result is proved:
\begin{tho}\label{EXCKN}{\sl For all $u_0\in J^2(\R^3)$ there exists a suitable weak solution.}\end{tho}
As a consequence of the inequality \eqref{SEI} and of the existence theorem one gets
\begin{coro}\label{EE} {\sl A suitable weak solution enjoys the strong energy inequality:
\be\label{EEII}\dm u(t)\dm_2^2+2\intll st \dm\nabla u(\tau)\dm_2^2d\tau\leq \dm u(s)\dm_2^2,\; \mbox{ for all } t\geq s, \mbox{ for } s=0 \mbox{ and a.e. in } s\geq 0\,.\ee Moreover for all $s$ such that \rf{EEII} holds we get\be\label{EEII-I}\lim_{t\to s^+}\dm u(t)-u(s)\dm_2=0\,.\ee}
\end{coro} 
Let us recall the definition of singular point for a weak solution. 
\begin{defi}\label{RP}{\sl We say
that $(t,x)$ is a singular point
for a weak solution $(u,\pi_u)$ if
$u\notin L^\infty$ in any
neighborhood of $(t,x)$; the
remaining points, where $u\in
L^\infty(I(t,x))$ for some
neighborhood  $I(t,x)$, are
called regular. }\end{defi}
\begin{defi}\label{RSD}{We say that $u$ is a  regular solution in $(t_0,t_1)\times\OO\subseteq (0,T)\times\R^3$  if $u$ is a weak solution, for some $q>1$, $u_t\in L^q_{\ell oc}((t_0,t_1)\times\OO))$   and, for all $\delta>0$, $u\in L^\infty((t_0+\delta,t_1-\delta)\times\OO)$\,.}\end{defi} It is known that a regular solution in $(t_0,t_1)\times\OO$ is smooth on   compact subsets contained  in $(t_0,t_1)\times\OO$, see e.g.  \cite{S}.
\par Following \cite{CKN} we introduce the parabolic cylinders 
\be\label{Qr}Q_r=Q_r(t,x):=\{(\tau,y):t-r^2<\tau<t\mbox{ and }|y-x|<r\},\ee and 
\be\label{Qrs}Q_r^*:=Q_r^*(t,x):=\{(\tau,y):t-\frac78r^2<\tau<t+\frac18r^2\mbox{ and }|y-x|<r\},\ee and, for $r\in (0,t^\frac12)$,  we set
\be\label{MCKN}M(r)=M(t,x,r):=r^{-2}\int\hskip-0.2cm\intl{Q_r}(|u|^3+|u||\pi_u|)dyd\tau+r^{-\frac{13}4}\intll{t-r^2}t\big(\hskip-0.15cm\intl{|x-y|<r}\hskip-0.2cm|\pi_u|dy\big)^\frac54d\tau\,,\ee
with $Q_r$  as in \eqref{Qr}.
\par In paper \cite{CKN}, in connection with the regularity of a suitable weak solution,  the authors furnish two regularity criteria.  The first is 
Proposition\,1 (or Corollary\,1, p.776) on p.775\,:
\begin{prop}\label{RCKN}{\sl Let $(u,\pi_u)$ be a suitable weak solution in some parabolic cylinder $Q_r(t,x)$. There exist $\vep_1>0$ and $c_0>0$ independent of $(u,\pi_u)$ such that, if \be\label{CKNSC}M(t,x,r)\leq \vep_1,\ee then
\be\label{RCKNI}|u(\tau,y)|\leq c_1^\frac12r^{-1},\mbox{ a.e. in }(\tau,y)\in Q_{\frac r2}(t,x),\ee where $c_1:=c_0\vep_1^{\frac23}$. In particular, a suitable weak solution $u$ is regular in $Q_{\frac r2}(t,x)$. }\end{prop} 
In \cite{CKN} this result is used to prove another regularity criterion, that is Proposition\,2 on p.776: 
\begin{prop}\label{RICKN}{\sl There is a constant $\vep_3>0$ with the following property. If $(u,\pi_u)$ is a suitable weak solution in some parabolic cylinder $Q_r^*(t,x)$ and 
$$\limsup_{r\to0}r^{-1}\int\hskip-0.2cm\intl{Q_r^*}
|\nabla u|^2dyd\tau\leq \vep_3\,,$$ then $(t,x)$ is a regular point.}\end{prop}
The above criterion is employed to get the following two main results (respectively, Theorem B on page 772 and Theorem D on page 774 in \cite{CKN}): 
\begin{tho}\label{BCKN}{\sl For any suitable weak solution the set $\mathbb S$ of singular points has one-dimensional parabolic Hausdorff measure equal to zero.}
\end{tho} 
\begin{tho}\label{DCKN}{\sl There exists an absolute constant $L_0>0$ with the following property. If $u_0\in J^2(\R^3)$, and if \be\label{CKN-L}\dm u_0|x|^{-\frac12}\dm_2=L<L_0\,,\ee then there exists a suitable weak solution to \rf{NS} which is regular in the region $$\{(t,x):|x|^2<t(L_0-L)\}\,.$$}\end{tho} 
 There is a difference in the meaning of the above theorems. Theorem\,\ref{BCKN} gives a geometric measure of the possible set S of singular points.    Theorem\,\ref{DCKN} it furnishes  the existence of a suitable weak solution to \rf{NS} having finite the following scaling invariant metric:
\be\label{WCKN}\sup_{0<\tau<t}\intl{\{\tau\}\times\R^3}|u|^2|x|^{-1}dx<\infty\,,\quad\intll0t\intl{\R^3}|\nabla u|^2|x|^{-1}dxd\tau<\infty\,\;t>0\,,\ee hence $x=0$ is   regular   for $t>0$. 
\par 
Finally,  as a corollary of the latter result, in \cite{CKN} the authors prove the following  (Corollary p. 820 in \cite{CKN}):
\begin{tho}\label{DCKN-C}{\sl Let $(u,\pi_u)$ be a suitable weak solution assuming initial data $u_0$. Suppose that $\dm\nabla u_0\dm_{L^2(|x|>R)}<\infty$. Then, there exists a $R_0>R$ such that, for all $\delta>0$, $u\in L^\infty((\delta,\infty)\times\{x:|x|>R_0\})$\,.}\end{tho} 
\subsection{\label{OT}\bf The aims of this note.} We work in the setting of the results of Theorem\,\ref{DCKN}   and Theorem\,\ref{TV} (below) already proved in \cite{CMRR}. Both these theorems work with a scaling invariant norm   that leads to   \rf{WCKN} provided that  at the initial instant the weighted norm, that is \rf{CKN-L}, \be\label{WN}\mathscr E(u_0,x):=\intl{\R^3}|u_0|^2|x-y|^{-1}dy\,,\;x\in\R^3\,,\ee is small in a suitable sense. The consequence of the smallness is the existence of a regular solution global in time.\par In this note we study the existence of a suitable weak solution that, at least locally in time, satisfies the regularity criterion of Proposition\,\ref{RCKN} and, as a consequence, is locally a regular solution. Also in this case the result follows from the assumption that the weighted norm \rf{WN} of the initial data is finite, but, contrary to  Theorem\,\ref{DCKN}   and Theorem\,\ref{TV}, we do not require smallness.  As a consequence we are able to deduce the regularity only locally in time.
\begin{tho}\label{TI}{\sl Let $u(t,x)$ be a suitable weak solution. Assume that for $x\in \mathbb E\subseteq\R^3$ there exists $v_0\in J^{1,2}(\R^3)$ such that\be\label{MA}\psi(x):=\intl{\R^3}\frac{|u_0(y)-u_0(y)|^2\hskip-0.1cm}{|x-y|}\hskip0.1cmdy<\frac1{(4c)^2\hskip-0.1cm}
\,,\ee where the constant $c$ is independent of $u_0,\,x,\,v_0$\,. Then there exist a $\delta\in[0,1)$ and $t>0$ such that \be\label{TI-I}u\in L^\infty(Q_{\big(\!\frac{(1-\delta)s}{4}\!\big)^\frac12}\mbox{$(\frac76s,x)$})\,\mbox{ for all }s\in(0,t)\,.\ee In particular,  if $(\tau,y)\in Q_{\big(\!\frac{(1-\delta)s}{4}\!\big)^\frac12}(\frac76s,x)$ is a Lebesgue point, then \be\label{TI-II}|u(\tau,y)|\leq \ov c\tau^{-\frac12}\,.\ee}
\end{tho} 
\begin{coro}\label{TIN}{\sl  Let $u(t,x)$ be a suitable weak solution. Then, for all $\sigma$ of validity of the weighted energy inequality {\rm(1.2)} there exists a set $\mathbb E\subseteq \R^3$, with $\mathbb \R^3-\mathbb E$ having zero Lebesgue measure, enjoying the property: for all $x\in\mathbb E(\sigma)$, there exist a $\delta\in [0,1)$ and $t>\sigma$ such that \be\label{CTI-I}u\in L^\infty(Q_{\big(\!\frac{(1-\delta)s}{4}\!\big)^\frac12}\mbox{$(\sigma+\frac76s,x)$})\,\mbox{ for all }s\in(0,t)\,.\ee In particular,  if $(\tau,y)\in Q_{\big(\!\frac{(1-\delta)s}{4}\!\big)^\frac12}(\sigma+\frac76s,x)$ is a Lebesgue point, then \be\label{CTI-II}|u(\tau,y)|\leq \ov c(\tau-\sigma)^{-\frac12}\,.\ee}\end{coro}
 We give some comments. 
\par Firstly we observe  that Theorem\,\ref{TI} 
 seems similar to Theorem\,\ref{DCKN}. The difference is in the fact that we do not require condition \rf{CKN-L} to the initial data, but the weaker condition \rf{MA}, that is almost everywhere satisfied by means of $u_0\in J^2(\OO)$. The theorem  establishes a result of local regularity for a suitable weak solution of \rf{NS}.   The local character is expressed in \rf{TI-I} either by the fact that the solution is $L^\infty$ just on the parabolic cylinder, and by the fact that the height of the cylinder depends on $x$, through $t(x)$.\par Estimate \rf{TI-II} (resp. \rf{CTI-II}) expresses in what way   the solution can be singular in $t=0$ (resp. in $\sigma$) provided that $x\in \mathbb E$ (resp. $x\in \mathbb E(\sigma)$).  
\par {In the  way specified below, the set $\mathbb E$ represents the new aspect of our result of local regularity stated with an initial data in $J^2(\R^3)$. Actually, if we consider $u_0\in J^2(\R^3)$, then the Riesz potential \be\label{RP}\mathscr E(u_0,x):=\intl{\R^3}\frac{u_0^2(y)}{|x-y|}dy\ee is well posed a.e. in $x\in\R^3$. This claim is consequence of the fact that, by the Hardy-Littlewood-Sobolev theorem, the following transformation is well defined:
\be\label{NORM}u_0^2\in L^1(\R^3)\to \mathscr E(u_0,x):=\intl{\R^3}\frac{u_0^2(y)}{|x-y|}dy\in L(3,\infty)(\R^3).\ee  Hence it is almost everywhere finite. Denoting by $\{u^k_0\}$ a sequence of smooth functions  converging to $u_0$ in $L^2(\R^3)$, for example the mollified of $u_0$, for $x\in\R^3$ and $k\in\N$ we define  the sequence
\be\label{PSI}\psi^k(x):= \intl{\R^3}\frac{|u_0-u_0^k|^2}{|x-y|}dy\,.\ee By Hardy-Littlewood-Sobolev theorem (see Lemma\,\ref{AP}), it is easy to verify that the sequence $\{\psi^k\}$ converges to zero almost everywhere in $x\in\mathbb E\subseteq\R^3$.    This makes satisfied almost everywhere in $x$ the assumption \rf{MA} and $\mathbb E$ is the set indicated in Corollary\,\ref{TIN}. We prove that for any $x\in \mathbb E$ there exists a $t(x)>0$ such that $M(\frac 76 s,x,r)\leq \varepsilon_1$ for suitable $r$ and for any $s\in (0,t(x))$. This result, by means of Proposition \ref{RCKN}, ensures the regularity in 
$Q_{\frac r2}(\frac 76 s, x)$, for any $s\in (0,t(x))$. Therefore, if we denote by $\mathbb S_x$ the projection onto $\R^3$ of the set $\mathbb S$ of singular points given in Theorem\,\ref{BCKN} (whose one-dimensional Hausdorff measure is zero from the same theorem), throughout Corollary\,\ref{TI} we can claim that $\mathbb S_x\subseteq \R^3\setminus\mathbb E$.   This last claim makes clear that we do not improve the regularity exhibited in \cite{CKN} (according with the result proved in \cite{scheffer2}), but we investigate on the existence of a possible size, as function of $x$ belonging to $\mathbb E$, of the parabolic neighborhood  of regularity of a weak solution. In Corollary\,\ref{TIN} it is claimed a dependence on $\sigma$ of the set $\mathbb E$: this is due to the fact that  we have to employ both \rf{SEI} and the right-continuity   in $L^2$-norm of the weak solution.}\par  The following results are two main consequences of Theorem\,\ref{TI}.
 \begin{tho}\label{TIIN}{\sl Let $u(t,x)$ be a suitable weak solution. Assume the existence of $\O\subseteq \R^3$ and  $v_0\in J^{1,2}(\R^3)$ such that \be\label{TII-I}\psi(x)<\frac 1{(4c)^2}\mbox{ uniformly in }x\in \O\,.\ee Then there exists a $T_0$ such that \rf{TI-I}, and \rf{TI-II}, hold for all $(s,x)\in (0,T_0)\times \O$. }\end{tho} 
 We observe that if $\OO\equiv\R^3$ then Theorem\,\ref{TIIN} gives the existence of a regular solution $(u,\pi_u)$ on $(0,T_0)\times\R^3$. 
 \begin{coro}\label{SevEgo}{\sl  Let $u(t,x)$ be a suitable weak solution. For any $B(R)$ and for any $\varepsilon>0$, there exists a set $\O_\varepsilon\subset B(R)$, with $meas(B(R)\setminus \O_\varepsilon)<\varepsilon$, and there exists a $T_0(\varepsilon)>0$ such that  \rf{TI-I} holds for all $(s,x)\in (0,T_0(\varepsilon))\times  \O_\varepsilon$.}
\end{coro} 
\begin{tho}\label{TV}{\sl Let $u(t,x)$ be a suitable weak solution, and assume also that  ${\underset x {\mbox{ess} \sup}}\,\mathscr E(u_0,x)$ is sufficiently small. Then,  $(u,\pi_u)$ is regular for all $t>0$ and it is unique up to a function $c(t)$ for the pressure field.}\end{tho}  
The last  theorems  are the regular solutions counterpart   of Theorem\,\ref{TI} and Corollary\,\ref{TIN}, provided that the assumptions on the data are stronger than the simple assumption $u_0\in J^2(\R^3)$. The theorems work in the light  of the scaling invariant weighted norm \rf{RP}.\par 
Theorem \ref{TIIN} establishes a local existence result stated by requiring a ``suitable closeness'',   in the weighted norm \rf{RP}, of the initial data $u_0\in L^2(\R^3)$ to a  smooth function $v_0$. As the existence is achieved on the element $v_0$ of the approximation  which is close to $u_0$ in the metric \rf{RP}, we are not able to give a size of $T_0$ by means of $u_0$, but $(0,T_0)$ is just  ({\it a priori}) a   subinterval of existence of the smooth solution $(v,\pi_{v})$ corresponding to  $v_0$. In this connection we point out that the above question on the size of $T_0$ is the same that we meet assuming the data $u_0$ in $J^3(\OO)$ or in $\mathbb L^3(\OO)\subset L(3,\infty)$, respectively completion of $\mathscr C_0(\OO)$ in $L^3(\OO)$ and in $L(3,\infty)(\OO)$. Both these spaces are scaling invariant and in order to prove the existence local in time we need an auxiliary function, say $\ov u_0,$ which is close to $u_0$ in the metric of $L^3$ or $L(3,\infty)$ and $\ov u_0\in X$, where $X$ is a function space adequate  to ensure  the existence of a regular solution on some interval $(0,T_0)$. This is an aspect   developed with details in \cite{MRM}. We conclude that in the statement of Theorem\,\ref{TIIN} we can substitute $J^{1,2}$ with any space $X$ which is suitable to ensure the existence of a regular solution corresponding to $v_0$. 
\par Corollary\,\ref{SevEgo} makes operational condition \eqref{TII-I} on a suitable subdomain. Indeed the existence of the domain $\OO_\vep\subseteq B(R)$ follows  from the construction of a sequence $\{\psi^k\}$  almost everywhere converging to zero and the Severini-Egorov theorem. 
\par  Theorem\,\ref{TV} furnishes a global existence result just requiring a smallness condition. It is also an immediate consequence of our previous result in \cite{CMRR}.
\section{Preliminaries}
Below we
recall 
some 
results which are fundamental for our
aims.  \begin{lemma}\label{WI}{\sl Suppose that
$|x|^\beta u\in L^2(\R^3)$ and
$|x|^\alpha \nabla u\in L^2(\R^3)$.
Also
\begin{itemize}\item[i)] $r\geq2$, $\gamma+\frac3r>0$, $\alpha+\frac32>0$, $\beta+\frac32>0$, and $a\in[\frac12,1]$,
\item[ii)]
$\gamma+\frac3r=a(\alpha+\frac12)+(1-a)(\beta+\frac32)$
(dimensional balance), \item[iii)]
$a(\alpha-1)+(1-a)\beta\leq\gamma\leq
a\alpha+(1-a)\beta$.
\end{itemize}
Then, with a constant $c$
independent of $u$,
 the following inequality holds:
\be\label{WII}\dm |x|^\gamma
u\dm_r\leq c\dm |x|^\alpha\nabla
u\dm_2^a\dm |x|^\beta
u\dm^{1-a}_2.\ee
 }\end{lemma}
\Pr See \cite{CKN}
Lemma\,7.1\,.\chiu
 \begin{lemma}\label{CZW}{\sl Assume that $\mathbb K$ is
a singular bounded transformation from $L^p$ into $L^p$,
$p\in(1,\infty)$,  of
Calder\'on-Zigmund kind. Then,
$\mathbb K$ is also a bounded
transformation from $L^p$ into $L^p$
with respect to the measure
$(\mu+|x|)^\alpha dx,\mu\geq 0,$ provided that
$\alpha\in
(-n,n(p-1))$.}\end{lemma}\Pr See
\cite{ST} Theorem\,1.\chiu
\begin{lemma}\label{RFP}{\sl Assume that $(u,\pi_u)$ is a suitable weak solution. Then the pressure field admits the following representation formula
\be\label{PI} \displ\pi_u(t,x)=-D_{x_i}D_{x_j}\intl{\R^3}\mathcal E(x-y)u^i(y)u^j(y)dy\,,\mbox{ a.e. in }(t,x)\in (0,\infty)\times\R^3\,,\ee and the following holds:\be\label{PIC} \pi_u(t,x)\in L^\frac53(0,T;L^\frac53(\R^3))\,. \ee}\end{lemma}
\Pr See
\cite{CMRR} Lemma\,2.4. Moreover, since $u^2\in L^\frac53(0,T;L^\frac53(\R^3))$ estimate \rf{PIC} easily follows.\chiu\begin{lemma}\label{RS}{\sl For all $v_0\in J^{1,2}(\R^3)$  there exists a unique regular solution $(v,\pi_v)$ to problem \rf{NS} on some interval $(0,T)$ such that
\be v\in C([0,T);J^{1,2}(\R^3)),\;v_t,\,D^2v,\,\nabla\pi_v\in L^2(0,T;L^2(\OO))\,,\ee where $T\geq c{\dm \nabla u_0\dm_2^{-4}}$.}\end{lemma}\Pr The result is due to Leray, see \cite{L}.\chiu   For $\mu\geq 0$ we define the functionals \be\label{FU}\mathscr E(v,t,x,\mu):=\intl{\R^3}\frac{|v(t,y)|^2}{(|x-y|^2+\mu^2)^\frac12 } dy\,,\quad \mathscr D(v,t,x,\mu):=\intl{\R^3}\frac{|\nabla v(t,y)|^2}{(|x-y|^2+\mu^2)^\frac12  } dy\,,\ee 
and set   $$p(y):=(|x-y|^2+\mu^2)^{-\frac12}.$$  When no confusion arises, we omit some or all the dependences on $(v,t,x,\mu)$. For  $\mu\geq0$, we call \be\label{WEI}  \mathscr E(v,t,x,\mu)+\intll0t \mathscr D(v,t,x,\mu)d\tau \ee {\it weighted energy}.\begin{lemma}\label{WERS}{Let $(v,\pi_v)$ be the regular the solution of Lemma\,\ref{RS}. Then, for all $\mu>0$, the following weighted energy relation and weighted energy inequality hold:
\be\label{WERSI}\ba{l}\displ \mathscr E(v,t,x,\mu)+ 2\!\intll0t \!\mathscr D(v,\tau,x,\mu)d\tau\!+ 3\mu^2\!\!\intll0t\!\intl{\R^3}\!\frac{v^2(\tau,y)}{(|x\!-\!y|^2+\!\mu^2)}\null_\frac52dy d\tau\!=\displ \mathscr E(v,0,x,\mu)\VS \hfill+\!\intll0t\!\intl{\R^3}v\otimes v\cdot v\otimes\nabla p dyd\tau+2\!\intll0t\!\intl{\R^3}\! \pi_v  v\cdot \nabla pdyd\tau,\ea\ee 
\be\label{Agg1}
\ba{l}\displ\mathscr E(t,x,\mu)+\intll0t\mathscr D(\tau,x,\mu)d\tau+3\mu^2\intll0t\intl{\R^3}\frac{v^2(\tau,y)}{(|x-y|^2+\mu^2)}\null_{\frac52}dyd\tau\VS\hskip4cm\leq \mathscr E(0,x,\mu)+c\intll0t\mathscr E(\tau,x,\mu)\dm\nabla v(\tau)\dm_2^4d\tau\,,\ea\ee 
for all $t\in\![0,T)$  and $x\in\!\R^3$. }\end{lemma}
\Pr Identity \eqref{WERSI} can be formally obtained 
by multiplying equation \rf{NS}$_1$ by $vp$ and  integrating by parts on $(0,t)\times\R^3$. Let us show that it is well posed for any $\mu> 0$. 
We start by remarking that in our hypotheses on $v_0$ we get $\mathscr E(0,x,\mu)<\infty$ for all $x\in\R^3$ and $\mu\geq0$.  By multiplying equation \rf{NS}$_1$ by $vp$ and  integrating by parts on $(0,t)\times\R^3$, we obtain \be\label{WERSII}\ba{l}\displ\mathscr E(t,x,\mu)+2\intll0t\mathscr D(\tau,x,\mu)d\tau+3\mu^2\intll0t\intl{\R^3}\frac{v^2(\tau,y)}{(|x-y|^2+\mu^2)}\null_\frac52 dyd\tau=\mathscr E(0,x,\mu)\VS\hskip 1cm -2\intll0t\!\intl{\R^3}(v\cdot \nabla v)\cdot vpdyd\tau-2\intll0t\!\intl{\R^3} \nabla\pi_v\cdot  v\cdot pdyd\tau\VS\hskip4cm=:\mathscr E(0,x,\mu)+2\!\intll0t\!(J_1\!+\!J_2)d\tau\,.\ea\ee
Let us show that the right-hand side is well defined. 
 Applying H\"older's inequality and inequality \rf{WII}, we get
$$|J_1|\leq \dm v(|x-y|^2+h^2)^{-\frac14}\dm_4^2\dm \nabla v\dm_2\leq \mathscr E^\frac14\mathscr D^\frac34\dm\nabla v\dm_2\leq \frac14\mathscr D+c\mathscr E\dm\nabla v\dm_2^4\,.$$
From the representation formula \rf{PI}, after integrating by parts, we get
$$ \nabla\pi_v(t,x)= \nabla\!\intl{\R^3} D_{y_j}\mathcal E(x-y)v^i(y)D_{y_i}v^j(y)dy\,.$$ Hence, applying H\"older's inequality and  employing Lemma\,\ref{CZW}, we deduce
$$\ba{ll}\displaystyle\vspace{1ex}
|J_2|\leq  \dm \nabla \pi_v(|x-y|^2+\mu^2)^{-\frac14}\dm_\frac43\dm v(|x-y|^2+\mu^2)^{-\frac14}\dm_4\\\hskip2cm \displaystyle\leq c\dm v\cdot\nabla v(|x-y|^2+\mu^2)^{-\frac14}\dm_\frac43\dm v(|x-y|^2+\mu^2)^{-\frac14}\dm_4\,.\ea$$ Applying again H\"older's inequality and subsequently \rf{WII}, we deduce the following estimate:
$$|J_2|\leq c \dm v(|x-y|^2+\mu^2)^{-\frac14}\dm_4^2\dm \nabla v\dm_2\leq c\mathscr E^\frac14 \mathscr D^\frac34\dm \nabla v\dm_2\leq \frac14\mathscr D+c\mathscr E\dm\nabla v\dm_2^4\,.$$ Hence from \rf{WERSII} and via estimates for terms $J_1$ and $J_2$ we obtain the integral inequality \eqref{Agg1}, 
from which, thanks to the regularity of $v$,  it is easy to deduce that \rf{WERSI} holds for all $\mu> 0$ and for all $t\in[0,T)$.  
\chiu
\begin{lemma}\label{AP}- {\sl Let
$u_0\in J^2(\R^3)$. There exists a set $\mathbb E$ such that $\R^3-\mathbb E$ has zero Lebesgue measure,    and for all $x\in\mathbb E$ and for all $\eta>0$   there
exists a $\ov u_0\in J^{1,2}(\R^3)$ such that 
\be\label{EXAPI}
 \intl{\R^3}\frac{|u_0-\ov u_0|^2\hskip-0.15cm\null}{|x-y|}\hskip0.15cmdy<\eta\,.\ee Moreover, for all $R>0$ and $\vep>0$ there exists $\OO_\vep\subseteq \mathbb E$  such that $meas(B(R)-\OO_\vep)<\vep$  and \be\label{EXAPI-I}
 \intl{\R^3}
\frac{|u_0-\ov u_0|^2\hskip-0.15cm\null}{|x-y|}\hskip0.15cmdy<\eta\,\mbox{ uniformly in }x\in \OO_\vep\,.\ee
}\end{lemma} \Pr We denote by $\{u^k_0\}$ the   mollified functions of $u_0$. It
is known that $\{u_0^k\}\subset
C^\infty(\R^3)\cap
J^{1,2}(\R^3),$ and
$\{u_0^k\}$ converges to
$u_0$ in $L^2$-norm. 
   For all $k\in\N$, we define \rf{PSI}, that is
 $$\psi^k(x):= \intl{\R^3}\frac{|u_0-u_0^k|^2\hskip-0.15cm}{|x-y|}\hskip0.15cmdy<\infty\,.$$ By the
Hardy-Littlewood-Sobolev theorem we
get, for $r\in[1,3)$, $$\dm
\psi^k\dm_{L^r(\mathbb K)}\leq
c(r,\mathbb K)\dm
u_0^k-u_0\dm_2^2,\mbox{ for
all compact set }\mathbb
K\subset\R^3.$$ Hence, the sequence
$\{\psi^k\}$ converges to zero in
$L^r(\mathbb K)$, for all $r\in
[1,3)$. In particular, there
exists a subsequence
$\{\psi^{k_j}\}$ which converges to
zero almost everywhere in $x\in
\mathbb  K$. We denote by $\{\mathbb K_\nu\}$ a sequence of compact sets such that $\mathbb K_\nu\subset \mathbb K_{\nu+1}$ and  ${\underset {\nu\in\N}\cup}\mathbb K_\nu=\R^3$. By virtue of the above convergence, we denote $\mathbb E_\nu\subseteq\mathbb K_\nu$ the set of  the convergence almost everywhere of the sequence $\{\psi^{k_j}\}$. Then, by means of  Cantor's diagonal method, we construct a sequence $\{\psi^\ell\}$ which converges to $0$ for all $x\in \mathbb E:={\underset {\nu\in\N}\cup}\mathbb E_\nu$. 
 Hence for all $x\in\mathbb E$ and $\eta>0$ there exists a $\psi^{\ov \ell}\in \{\psi^{\ell}\}$ such that $\ov u_0:=\ov u_0^\ell$ verifies \rf{EXAPI}. Property \rf{EXAPI-I} is a consequence of the above construction and of the Severino-Egorov theorem. 
The lemma is completely proved.\chiu
\section{Local in time weighted energy inequality
for a suitable weak solution}
In this section we prove that any suitable weak solution admits at least locally in time a weighted energy inequality with $\mu=0$. Actually, the following lemma holds
\begin{lemma}\label{WERLT}{\sl Let $(u,\pi_u)$ be a suitable weak solution.  Let $x$, $v_0$ and $c$ as in Theorem\,\ref{TI}. Then 
there exists a $t^*(x)>0$ such that
\be\label{WEILI}
\mathscr E(u,t,x)+\frac12\intll0t\mathscr D(u,\tau,x)d\tau\leq N<\infty,\mbox{ for all }t\in[0,t^*(x)) \,,\ee
with $\mathscr E(u,t,x)$ and $\mathscr D(u,\tau,x)$ defined in \eqref{FU}.
}\end{lemma}
\Pr The proof  of estimate \rf{WEILI} reproduces in a suitable way an idea employed in \cite{CMSD}. This idea follows the Leray-Serrin arguments employed for the proof of the energy inequality in strong form. The proof  is achieved by means of five steps. We set $w:=u-v$ and $\pi_w:=\pi_u-\pi_v$, where $(u,\pi_u)$ is the suitable weak solution and $(v,\pi_v)$ the regular solution corresponding to $v_0$ and furnished by Lemma\,\ref{RS}. The first four steps are devoted to prove the following inequality
\be\label{WERLTI}
\mathscr E(w,t,x,\mu)+\frac12\intll0t\mathscr D(w,\tau,x,\mu)d\tau\leq \frac{1}{8c^2},\mbox{ for all }t\in[0,t^*(x))\mbox{ and }\mu>0\,.\ee
\par{\it Step 1.} - 
 We start proving that for all $t>0$
\be\label{DIWE}\ba{l}\displ\mathscr E(t,x,\mu)
+2\intll0t\intl{\R^3}\mathscr D(\tau,x,\mu)d\tau+\!3\mu^2\!\!\!\intll0t \!\!\intl{\R^3}\!
\frac{|u(\tau)|^2}{(|x\!-\!y|^2\!+\!\mu^2)}\null_{\!\frac52}dyd\tau
\VS\hskip1cm\leq
\mathscr E(0,x,\mu)
+\!\intll0t\!\!\intl{\R^3}\!\frac{|u(\tau)|^2\,u\cdot(x\!-\!y)}{(|x\!-\!y|^2+\mu^2)^\frac 32\hskip-0.5cm\null}\,
dyd\tau \!+2\!\!\intll0t\!\!\intl{\R^3}\!
\frac{\pi_{u}(\tau)u(\tau)\!\cdot\!(x\!-\!y)}{(|x\!-\!y|^2+\mu^2)^\frac 32\hskip-0.5cm\null}\,dyd\tau.\ea\ee
  In the energy inequality \rf{SEI} we set $\phi(\tau,y):=(|x-y|^2+\mu^2)^{-\frac12}h_R(y)k(\tau)\in C_0^\infty(\R\times\R^3)$, with $h_R$ and $k$    such that
$$h_R(y):=\left\{\ba{ll}1&\mbox{if }|y|\leq R\\\in (0,1)&\mbox{if }|y|\in (R,2R)\\0&\mbox{for }|y|\geq 2R\,,
\ea\right.\mbox{ and }k(\tau):=\left\{\ba{ll}1&\mbox{if }|\tau|\leq t\\\in (0,1)&\mbox{if }|\tau|\in (t,2t)\\0&\mbox{for }|\tau|\geq 2t\,.
\ea\right.$$ We get
\be\label{WEIIIA}\ba{l}\displ\!\! \intl{\R^3}\!\frac{|u(t)|^2h_R}{(|x\!-\!y|^2\!+\!\mu^2)}\null_{\!\frac12}dy
+\!2\!\!\intll0t\!\!\intl{\R^3}\!\frac{|\nabla u(\tau)|^2h_R}{(|x\!-\!y|^2\!+\!\mu^2)}\null_{\!\frac12}dyd\tau\!+\!3\mu^2\!\!\!\intll0t \!\!\intl{\R^3}\!
\frac{|u(\tau)|^2h_R}{(|x\!-\!y|^2\!+\!\mu^2)}\null_{\!\frac52}dyd\tau
\VS\hskip1cm\leq
\intl{\R^3}\!\frac{|u_0|^2h_R}{(|x-y|^2\!+\!\mu^2)}\null_{\!\frac12}dy
+\!\intll0t\!\intl{\R^3}\!\frac{|u(\tau)|^2h_R\,u\cdot(x-y)}{(|x-y|^2+\mu^2)^\frac 32}\,
dyd\tau\VS\hskip1.5cm+2\intll0t\!\intl{\R^3}\!
\frac{\pi_{u}(\tau)h_Ru(\tau)\cdot\!(x-y)}{(|x-y|^2+\mu^2)^\frac 32}\,dyd\tau+F(t,R)\VS\hskip 4cm:=\intl{\R^3}\!\frac{|u_0|^2h_R}{(|x-y|^2\!+\!\mu^2)}\null_{\!\frac12}dy+I_1(t,x)+I_2(t,x)+F(t,R),\ea\ee where
$$\ba{l}\displ F(t,R):=\intll0t \intl{\R^3}|u|^2\Big[2\nabla h_R\cdot\nabla (|x-y|^2+\mu^2)^{-\frac12}+\frac{\Delta h_R}{(|x-y|^2+\mu^2)}\null_{\!\frac12}+\frac{u\cdot \nabla h_R}{ (|x-y|^2+\mu^2)}\null_{\frac12}\Big]dyd\tau\VS\hskip3cm+\intll0t \intl{\R^3}\frac{\pi_uu\cdot\nabla h_R}{(|x-y|^2+\mu^2)}\null_{\!\frac12}dyd\tau \,.\ea$$ Since $\pi_u,\,u^2\in L^\frac53(0,T;L^\frac53(\R^3))$, applying H\"older's inequality and employing the decay of $\nabla h_R,\,\Delta h_R$, for all $t>0$, we get   $F(t,R)=o(R)$. We estimate the terms $I_i,\,i=1,2$. Since $\mu>0$, by virtue of the integrability properties of a suitable weak solution, applying Lemma\,\ref{WI} we get
$$|I_1(t,x)|\leq 
  \intll0t\dm \frac{u}{(|x-y|^2+\mu^2)}\null_{\!\frac13}\dm_3 ^3d\tau \leq c\intll0t
\dm \frac{u}{(|x-y|^2+\mu^2)}\null_{\frac14}\dm_2
\dm \frac{\nabla u}{(|x-y|^2+\mu^2)}\null_{\frac14}\dm_2^2d\tau.$$  
For $I_2$, applying the H\"older's inequality and 
Lemma\,\ref{CZW}, we obtain
$$|I_2(t,x)|\leq c\intll0t\dm \frac{u}{(|x-y|^2+\mu^2)}\null_{\frac13}\dm_3
\dm \frac{\pi_{u}}{(|x-y|^2+\mu^2)}\null_{\frac23}\dm_\frac32d\tau\leq c\intll0t \dm \frac{u}{(|x-y|^2+\mu^2)}\null_{\frac13}\dm_3^3d\tau.$$ Hence, as in the previous case, applying  Lemma\,\ref{WI}, we get $$|I_2(t,x)|\leq c\intll0t
\dm \frac{u}{(|x-y|^2+\mu^2)}\null_{\frac14}\dm_2
\dm \frac{\nabla u}{(|x-y|^2+\mu^2)}\null_{\frac14}\dm_2^2d\tau.$$
Employing the estimates obtained for $I_i,i=1,2$, via the Lebesgue dominated convergence theorem, in the limit as $R\to\infty$, for all $t>0$ we deduce the  inequality \rf{DIWE}.  \vskip0.1cm\par {\it Step 2.} - In this step we derive a sort of Green's identity between solutions $(u,\pi_u)$ and $(v,\pi_v)$, where $(v,\pi_v)$ is the regular solution given in Lemma\,\ref{RS}, corresponding to the initial data $v_0\in J^{1,2}(\R^3)$. In the following $(0,T)$ is the interval of existence of $(v,\pi_v)$. We also recall that the regular solution $(v,\pi_v)$ is smooth for $t>0$.  We denote by $\lambda(\tau)$ a smooth cutoff function such that $\lambda(\tau)=1$ for $\tau\in[s,t]$ and $\lambda(\tau)=0$ for $\tau\in[0,\frac s2]$.\par
For all $t,s\in(0,T)$, we consider the weak formulation iii) of Definition\,\ref{WS} written with $\varphi=\lambda vp$: 
\be\label{WERLTII}\ba{l}\displ\intll
st\Big[(pu,v_\tau)-(p\nabla
u,\nabla
v)+(pu\cdot\nabla v,u)+(\pi_u,v\cdot\nabla
p)\Big]d\tau+(pu(s),v(s)
)\VS\hskip3cm =(pu(t),v(t)) +
\displ\intll
st\Big[ (\nabla
u,v\otimes\nabla
p)+(u\otimes u, v\otimes\nabla p) \Big]d\tau\, .
\ea\ee
We multiply equation \rf{NS}$_1$ written for $(v,\pi_v)$ by $u p$. After  integrating by parts on $(s,t)\times\R^3$, we get 
\be\label{WERLTIII}\ba{l}\displ\intll
st\Big[(pu,v_\tau)+(p\nabla
u,\nabla
v)+(pv\cdot\nabla v,u)-(\pi_v,u\cdot\nabla
p)\Big]d\tau \VS\hskip 5cm= 
\displ\intll
st\Big[(\nabla u
,v\otimes\nabla
p)+(u\cdot v,\Delta p)\Big]d\tau\, .
\ea\ee making the difference between formulas \rf{WERLTII} and \rf{WERLTIII} we get
$$\ba{l}\displ\intll
st\Big[ -2(p\nabla
u,\nabla
v)+(pu\cdot\nabla v,u)-(pv\cdot\nabla v,u)+(\pi_u,v\cdot\nabla
p)+(\pi_v,u\cdot\nabla
p)\Big]d\tau\VS\hskip1.6cm 
 =(pu(t),v(t)) -(pu(s),v(s))+\!
\displ\intll
st\!\Big[ (u\otimes u, v\otimes\nabla p)-(u\cdot v,\Delta p)\Big] d\tau\, ,
\ea$$ Since in a suitable neighborhood of $0$ all the terms of the last integral equation are continuous on the right, letting $s\to0^+$, we get
\be\label{WERLTIV}\ba{l}\displ\intll
0t\Big[ -2(p\nabla
u,\nabla
v)+(pu\cdot\nabla v,u)-(pv\cdot\nabla v,u)+(\pi_u,v\cdot\nabla
p)+(\pi_v,u\cdot\nabla
p)\Big]d\tau\VS\hskip1.6cm 
 =(pu(t),v(t)) -(pu(0),v(0))+\!
\displ\intll
0t\!\Big[ (u\otimes u, v\otimes\nabla p)-(u\cdot v,\Delta p)\Big] d\tau\, ,
\ea\ee
which furnishes the wanted Green's identity.\vskip0.1cm\par {\it Step 3.} - Setting $w:=u-v$ and $\pi_w:=\pi_u-\pi_v$, let us derive the following estimate \be\label{WERLTVI}\hskip-0.3cm\ba{l} \displ  \mathscr  E(w,t,x,\mu)\!+\!\!\intll0t\!\!\mathscr D(w,\tau,x,\mu)d\tau \leq\! \mathscr  E(w,0,x,\mu)\!    +c\!\!\intll0t\!\!\mathscr E^\frac12(w,\tau,x,\mu)\mathscr D(w,\tau,x,\mu)d\tau\VS\hskip5cm+ H(v,t,x,\mu),\mbox{\,for all }t\!\in\![0,T),x\!\in\!\R^3,\mu\!>\!0,\ea\ee with $$H(v,t,x,\mu):=c\!\intll0t\!\! \dm\nabla v(\tau)\dm_2^4d\tau\!+\!c\!\intll0t\!\!\mathscr E(v,\tau,x,\mu)\mathscr D(\tau,v,x,\mu)d\tau.$$  
  We remark that from the representation formula \rf{PI} and regularity of $v$ we get that \be\label{PI-PII}\ba{c}\pi_w=\pi^1+\pi^2\,,\VS\pi^1\!\!:= D_{x_j}\!\!\intl{\R^3}\!\!D_{y_i}\mathcal E(x\!-\!y)w^i(y)w^j(y)dy\mbox{\, and \,}\pi^2\!\!:=2  \!\intl{\R^3}\!\!D_{y_j} \mathcal E(x\!-\!y)w(y)\!\cdot\!\nabla v^j(y)dy.\ea\ee
We sum estimates \rf{WERSI} and \rf{DIWE}, then  we add twice formula \rf{WERLTIV}.  written for $s=0$. Recalling the definition of $(w,\pi_w)$ and formula \rf{PI-PII}, after a straightforward computation  we get
\be\label{WERLTV}\ba{c}\displ  \mathscr  E(w,t,x,\mu)+2\intll0t\!\mathscr D(w,\tau,x,\mu)d\tau+3\mu^2\intll0t\!\intl{\R^3}\!\frac{w^2(\tau,x)}{(|x-y|^2+\mu^2)}\null_\frac52d\tau\VS\leq \mathscr  E(w,0,x,\mu)+F_1(w,t,x,\mu)+F_2(w,v,t,x,\mu),\ea\ee
where   $$\ba{c}\displ F_1:=F_1(w,t,x,\mu):=\intll0t\!(w\otimes w,w\otimes\nabla p)d\tau+2\!\intll0t\!( \pi_1,  w\cdot \nabla p)d\tau\VS F_2:=F_2(w,v,t,x,\mu):= 2\intll0t(\pi_2,w\cdot\nabla p)d\tau-2\intll0t(w\cdot\nabla v,wp)d\tau+ \intll0t(v\cdot\nabla p,w^2) \,.\ea$$ The term $F_1$ admits the same estimate as $I_1$ and $I_2$ given in {\it Step 1}, hence we get
$$|F_1|\leq c \intll0t \mathscr E^\frac12(\tau,w,x,\mu)\mathscr D(\tau,w,x,\mu)d\tau\,\mbox{ for all }t\in(0,T),\,x\in\R^3,\,\mu>0.$$
For term $F_2$ we estimate the first two terms in a different way from the last. Taking  the representation formula of $\pi_2$ into account, we get
$$|\intll0t(\pi_2,w\cdot\nabla p)d\tau-2\intll0t(w\cdot\nabla v,wp)d\tau|=|\intll0t p\nabla\pi_2\cdot wdyd\tau+2\intll0t(w\cdot\nabla v,wp)dyd\tau |\,.$$ Hence, applying the same arguments employed in Lemma\,\ref{WERS} to estimate $J_1$ and $J_2$, we get
$$\ba{l}\displ|\intll0t(\pi_2,w\cdot\nabla p)d\tau-2\intll0t(w\cdot\nabla v,wp)d\tau|\leq 
\int_0^t\dm wp^\frac12\dm_4^2\dm \nabla v\dm_2d\tau\VS\hskip 0.3cm
  \leq \intll0t \mathscr E^\frac13(w,\tau,x,\mu)\mathscr  D(w,\tau,x,\mu)d\tau +c\intll0t \dm\nabla v(\tau)\dm_2^4d\tau\,, \mbox{ for all }t\in[0,T),\,x\in\R^3,\,\mu>0.\ea$$
For the last term in $F_2$, applying H\"older's inequality, we get
$$|\intll0t(v\cdot\nabla p,w^2)d\tau|\leq \int_0^t\dm wp^\frac12\dm_4^2\dm vp\dm_2^2d\tau\,.$$ By virtue of estimate \rf{WII}, applying Young's inequality we deduce:
$$\ba{ll}|\displ\intll0t(v\cdot\nabla p, w^2)d\tau|&\hskip-0.2cm\leq c\displ\intll0t\mathscr E^\frac14(w,\tau,x,\mu)\mathscr D^\frac34(w,\tau,x,\mu)\mathscr E^\frac14(v,\tau,x,\mu)\mathscr D^\frac14(v,\tau,x,\mu)d\tau\VSE\hskip-0.2cm\leq  \intll0t\mathscr E^\frac13(w,\tau,x,\mu)\mathscr D(w,\tau,x,\mu)d\tau+c\intll0t\mathscr E(v,\tau,x,\mu)\mathscr D(v,\tau,x,\mu)d\tau.\ea$$ Hence, we obtain  
$$\ba{l}\displ|F_2|\leq    2 \intll0t\mathscr E^\frac13(w,\tau,x,\mu)\mathscr D(w,\tau,x,\mu)d\tau+c\intll0t \dm\nabla v(\tau)\dm_2^4d\tau\VS\hskip3cm+c\intll0t\mathscr E(v,\tau,x,\mu)\mathscr D(v,\tau,x,\mu)d\tau, \mbox{ for all }t\in[0,T),\,x\in\R^3,\,\mu>0.\ea $$ Finally, applying Young's inequality, we get
$$\ba{l}\displ|F_2|\leq  \intll0t\mathscr D(w,\tau,x,\mu)d\tau+  c \intll0t\mathscr E^\frac12(w,\tau,x,\mu)\mathscr D(w,\tau,x,\mu)d\tau+c\intll0t \dm\nabla v(\tau)\dm_2^4d\tau\VS\hskip3cm+c\intll0t\mathscr E(v,\tau,x,\mu)\mathscr D(v,\tau,x,\mu)d\tau, \mbox{ for all }t\in[0,T),\,x\in\R^3,\,\mu>0.\ea $$
Estimates for $F_1,F_2$ and \rf{WERLTV}   furnish the integral inequality \rf{WERLTVI}.\vskip0.1cm
\par{\it Step 4.} - Deduction of estimate \rf{WERLTI}. \par
Under our assumptions on $x$, $v_0$ and $c$,  we have, a fortiori,
\be\label{HP}\mathscr E(w,0,x,\mu)<\frac{1}{(4c)}\null_2,\ \mbox{for all } \mu>0.\ee Moreover
by virtue of the regularity of the solution $(v,\pi_v)$ (see Lemma\,\ref{RS} and Lemma\,\ref{WERS}) there exists a $ t^*$ such that  
\be\label{HPN}H( t^*)<\frac1{(4c)}\null_2\,,\mbox{for all } \mu>0.\ee
 Let us deduce \rf{WERLTI} that for convenience of the reader we rewrite: \be\label{WEW}  \displ  \mathscr  E(w,t,x,\mu)\!+\mbox{$\frac1{2}$}\!\intll0t\!\!\mathscr D(w,\tau,x,\mu)d\tau <\mbox{$\frac1{8c^2}$}\,, \mbox{ for all }t\in [0, t^*), \,\mu>0. \ee  
Since $w=u-v$ is right continuous in $L^2$-norm in $t=0$, for all $\mu>0$ the same continuity  property holds for $\mathscr E(w,t,x,\mu)$. 
Therefore there exists 
a $\delta=\delta(\mu)>0$ such that \be\label{SLM}\mathscr E(w,t,x,\mu)<\frac 1{8c^2} \,, \mbox{ for all }t\in[0,\delta).\ee 
Hence the validity of estimates \rf{WERLTVI} and  \eqref{HP}-\eqref{HPN} yields  for any  $t \in [0,\delta)$
$$\mathscr  E(w,t,x,\mu)\!+\!\!\intll0{t}\!\!\mathscr D(w,\tau,x,\mu)d\tau <\!\frac 1{8c^2}+c\!\intll0{ t}\!\!\mathscr E^\frac 12(w,\tau,x,\mu)\mathscr D(w,\tau,x,\mu)d\tau,$$ that, thanks to \eqref{SLM}, gives \eqref{WEW} on $[0,\delta)$. \par Let us show that estimate \rf{SLM} holds for $t\in[0,t^*)$.
For all $\mu>0$, the function 
$$f(t,\mu):=\mathscr  E(w,0,x,\mu)\!    +c\!\intll0t\!\!\mathscr E^\frac12(w,\tau,x,\mu)\mathscr D(w,\tau,x,\mu)d\tau + H(v,t,x,\mu)$$
is uniformly continuous on $[0,t^*]$. Hence there exists $\eta=\eta(\mu)>0$ such that $$|t_1-t_2|<\eta\Rightarrow |f(t_1)-f(t_2)|<\frac1{8c}\null_2-\mathscr E(w,0,x,\mu)-H(t^*(x))\,.$$
We state that estimate \rf{SLM} and, consequently, estimate \eqref{WEW}, also holds for $t\in[\delta,\delta+\eta)$. Assuming the contrary, there exists $\ov t\in[\delta,\delta+\eta)$ such that 
\be\label{SLMI} \mathscr E(w,\ov t,x,\mu)>\frac 1{8c^2}\,.\ee On the other hand, the validity of \rf{WERLTVI} yields 
$$\mathscr  E(w,\ov t,x,\mu)\!+\!\!\intll0{\ov t}\!\!\mathscr D(w,\tau,x,\mu)d\tau \leq (f(\ov t)-f(\delta))+f(\delta)<\!\frac 1{8c^2}+c\!\intll0{\delta}\!\!\mathscr E^\frac 12(w,\tau,x,\mu)\mathscr D(w,\tau,x,\mu)d\tau .$$
Estimate \rf{SLM} allows to deduce that $$c\intll0{\delta}\!\!\mathscr E^\frac 12(w,\tau,x,\mu)\mathscr D(w,\tau,x,\mu)d\tau< \frac 1{\sqrt 8}\intll0{\ov t}\mathscr D(w,\tau,x,\mu)d\tau.$$ Hence the last two estimates imply $$\mathscr E(w,\ov t,x,\mu)<\frac 1{8c^2},$$ which is in   contradiction with \rf{SLMI}. Since the arguments are independent of $\delta$, the result holds for any $t\in [0, t^*(x))$, which proves \rf{WEW}\par{\it Step 5.} - Since $u=w+v$, via estimate \rf{Agg1} and  via estimate \eqref{WEW} we deduce, with obvious meaning of $N$ and $t^*(x)$ independent of $\mu$, the following inequality     $$ \intl{\R^3}\frac{|u(t,y)|^2}{(|x-y|^2+ \mu^2)}\null_\frac12+(1-\mbox{$\frac 1{\sqrt8}$})\intll0t\intl{\R^3}\frac{|\nabla u(t,y)|^2}{(|x-y|^2+\mu^2)}\null_\frac12dyd\tau\leq N\,,\mbox{ for all }t\in[0,t^*(x))\,.$$ The thesis is an easy consequence of estimate \eqref{WERLTI}  and the following remark:
the families of functions 
$$\big\{\intll0t\intl{\R^3}\frac{|\nabla u(t,y)|^2}{(|x-y|^2+\mu^2)}\null_\frac12dyd\tau\big\}\mbox{ and }\big\{\intl{\R^3}\frac{|u(t,y)|^2}{(|x-y|^2+\mu^2)}\null_\frac12dy\big\}$$ are monotone in $\mu>0$. Hence,   by virtue of the Beppo Levi's theorem, in the limit as $\mu\to0$, we deduce \rf{WEILI}.\chiu
\begin{coro}\label{WEILN}{\sl Let $(u,\pi_u)$ be a suitable weak solution. Let $\sigma\geq0$ such that \rf{SEI} is verified. Then there exists a set $\mathbb E\subseteq\R^3$, with $\R^3-\mathbb E$ having zero Lebesgue    measure, enjoying the property: for all $x\in\mathbb E(\sigma)$ there exists a $t^*(x)>0$ such that  
\be\label{WEILI-I}
\mathscr E(u,t,x)+(\mbox{$1-\frac1{\sqrt 8}$})\intll\sigma t\mathscr D(u,\tau,x)d\tau\leq N<\infty,\mbox{ for all }t\in[\sigma,\sigma+t^*(x)) \,.\ee}\end{coro}
\Pr    For all $\sigma\geq0$ for which $u$ verifies \rf{SEI}, via Lemma\,\ref{AP},  there exists a set $\mathbb E$ such that for $x\in \mathbb E$ and $\vep>0$ there exists a function $\ov u(\sigma)\in J^{1,2}(\R^3)$  that allows us to verify \eqref{MA} of Theorem\,\ref{TI}   with $u(\sigma)-\ov u(\sigma)$. As the assumptions of Lemma \,\ref{WERLT} are satisfied, the result follows. \chiu 
\section{Proof of Theorems\,\ref{TI}-\ref{TIIN} and Corollaries\,\ref{TIN}-\ref{SevEgo}.
}
To prove Theorem\,\ref{TI} we employ the result of Proposition\,\ref{RCKN}. To this aim, in the following Lemma\,\ref{RLI} we prove that, for a suitable $r>0$, estimate \rf{WEILI} of Lemma \,\ref{WERLT} implies condition \rf{CKNSC} of Proposition\,\ref{RCKN}. 
\begin{lemma}\label{RLI}{\sl    Let  the assumption of Lemma\,\ref{WERLT} be satisfied. Then,     there exists  $\delta>0 $ such that
\be\label{RLII}M(t,x,r)\leq \vep_1 \,,\mbox{ for all }r\in(0,[(1-\delta)t]^\frac12)\mbox{ and }t\in(0,t^*(x)).\ee 
with $t^*(x)$ given in Lemma\,\ref{WERLT}. 
 }\end{lemma}
\Pr  By virtue of our assumption, and by virtue of representation formula \rf{PI} and Lemma \ref{CZW},  a.e. in $t\in(0,t^*(x))$, we   get that 
\be\label{PIII} \dm \pi_u(t) |x-y|^{-\frac43}\dm_{\frac32}\leq c
\dm |u(t)| |x-y|^{-\frac23}\dm_{3}^2\,.\ee
Applying H\"older's inequality, from \rf{PIII} and from Lemma\,\ref{WI},  for all $t\in (0,t^*(x))$ and $t-r^2>0$, we have 
\be\label{PIV}\!\ba{ll}\displ r^{-2}\!\!\!\intll{t-r^2}t\,\intl{|x-y|<r}\!\!\!\!\Big[|u|^3+|v||\pi_u|\Big]dyd\tau\hskip-0.3cm&\displ\leq\! c\!\intll{t-r^2}t\!\!\Big[\!
\dm\frac{ u(\tau)}{ |x-y|}\null_{_{\!\!\frac23}}\!\dm_{3}^3+
\dm \frac{u(\tau)}{ |x-y|}\null_{_{\!\!\frac23}}\!\dm_{3}
\dm \frac{\pi_u(\tau) }{|x-y|}\null_{_{\!\!\frac43}}\!\dm_{\frac32}\!\Big]d\tau\VSE\leq \!c\!\intll{t-r^2}t\!\dm \frac{u(\tau)}{|x-y|
}\null_{\null_{\frac12}}\!\dm_2 \dm\frac{\nabla u(\tau)}{|x-y|}\null_{\null_{\!\frac12}}\!\dm_2^2d\tau\VSE=c\intll{t-r^2}t\mathscr E(\tau,x)^\frac12\mathscr D(\tau,x)d\tau=:N_1  .
\ea\ee 
Considering the second term on the right-hand side of $M(t,x,r)$ in \eqref{MCKN}, applying twice H\"older's inequality, \rf{PIII}, or all $t\in (0,t^*(x))$ and $t-r^2>0$, we get 
\be\label{PV} \ba {ll}\displ r^{-\frac{13}4}\hskip-0.3cm\intll{t-r^2}t\hskip-0.2cm\Big[\intl{|x-y|<r} \hskip-0.4cm|\pi_u(\tau,y)|dy\Big]^\frac54\!\!d\tau\hskip-0.3cm&\displ\leq 
cr^{-\frac 13}\hskip-0.3cm\intll{t-r^2}t\hskip-0.2cm
\Big[\dm \frac{\pi_u(\tau) }{|x-y|}\null_{_{\!\!\frac43}}\!\!\dm_{\frac32}\Big]^\frac54\!\!d\tau\leq cr^{-\frac13}\hskip-0.3cm\intll{t-r^2}t
\!\!\!\!\dm \frac{u(\tau)}{|x-y|
}\null_{\null_{\frac12}}\!\dm_2^\frac56 \dm\frac{\nabla u(\tau)}{|x-y|}\null_{\null_{\!\frac12}}\!\dm_2^\frac53\!d\tau\VSE\leq c\Big[\intll{t-r^2}t\mathscr E^\frac12(\tau,x)\mathscr D(\tau,x)d\tau\Big]^\frac56 =:N_2 .\ea\ee
Hence  \rf{PIV} and \rf{PV} imply  that
$$M(t,x,r)\leq N_1+N_2.$$ 
Employing estimate \rf{WEILI}, we get
$$N_1+N_2\leq cN^\frac12\!\!\!\intll{t-r^2}t\!\!\mathscr D(\tau,x)d\tau+
 \Big[cN^\frac12\!\!\!\intll{t-r^2}t\!\!\mathscr D(\tau,x)d\tau\Big]^\frac56\!, \mbox{ for all }t\!\in\!(0,t^*(x))\mbox{ and }t-r^2>0.$$ On the other hand   the function$$\intll t{t^*(x)}\mathscr D(\tau)d\tau\mbox{ is uniformly continuous on }[0,t^*(x)].$$ Hence there exists a   $\ov\delta\in(0,1)$   such that  
$$\Big[cN^\frac12\!\!\!\intll{(1-\ov\delta)t}t\!\!\mathscr D(\tau,x)d\tau\Big]^\frac56\!+cN^\frac12\!\!\!\intll{(1-\ov\delta)t}t\!\!\mathscr D(\tau,x)d\tau<\vep_1 \,\forall t\in (0,t^*(x)).$$
Hence the lemma is proved setting $\delta=1-\ov\delta$.\chiu 
Now we are in a position to prove the results of Theorem\,\ref{TI} and Theorem\,\ref{TIIN}. \vskip0.2cm
\par {\it Proof of Theorem\,\ref{TI}} - By virtue of Lemma\,\ref{WERLT}, for any $x$ satisfying the assumptions, estimate \rf{WEILI} holds on some interval $[0,t^*(x))$. 
Set $t(x):=\frac67 t^*(x)$, 
by virtue of  Lemma\,\ref{RLI}, 
there exists a $\delta>0$ such that $M(\frac 76s, x, r)\leq \varepsilon_1$, for all
$r\in(0,[(1-\delta)\frac 76s]^\frac12)$, $s\in(0,t(x))$. This, via Proposition\,\ref{RCKN}, implies the local regularity \rf{TI-I}, provided that  
$\delta\in (0,\frac 17)$\footnote{This condition ensures that we can choose $r=\sqrt s$, being  $(1-\delta)\frac 76s>s$.}. 
Finally, in order to prove \rf{TI-II} it is enough to observe that the point  $(s,x)$ belongs to $Q_{\left(\frac s4\right)^\frac12}(\frac76s,x)$ and, if $(s,x)$ is a Lebesgue point, then, via estimate \rf{RCKNI}, we can state \rf{TI-II}. The theorem is completely proved.\chiu
\vskip0.2cm
\par {\it Proof of Corollary\,\ref{TIN}} - By virtue of Corollary\,\ref{WEILN}, there exists a set  $\mathbb E(\sigma)$ such that for all $x\in\mathbb E(\sigma)$ estimate \rf{WEILI} holds on some interval $[\sigma,\sigma+t^*(x))$. Then one can conclude as in the proof of Theorem\,\ref{TI}. \chiu
\vskip0.2cm
\par{\it Proof of Theorem\,\ref{TIIN}} - Under the assumption of the theorem, Lemma\,\ref{WERLT} holds for any $x$ in $\O$, with $t^*(x)$ uniform in $\O$. The last claim is a consequence of the fact that in the definition of $\psi$ the smooth function  $v_0$ is independent of $x\in\OO$. Hence under our assumption \rf{TII-I} we have that both \rf{HP} and \rf{HPN} are uniform with respect to $x$. 
Setting $T_0:=t^*$, we write \rf{WEILI}   for $t\in[0,T_0)$ for all $x\in \Omega$. As a consequence, all the arguments employed for the proof of Theorem\,\ref{TI} work independently of $x\in \Omega$. The theorem is proved.\chiu\vskip0.2cm
\par{\it Proof of Corollary\,\ref{SevEgo}} - Fixed the ball $B(R)$ and given $\vep>0$, we can employ Lemma\,\ref{AP} which furnishes property \rf{EXAPI-I}. Hence the assumption of  Theorem\,\ref{TIIN} holds for any $x$ in $\O_\vep$. \chiu
 \vskip0.2cm 
 {\bf Acknowledgments} - {\small The research is performed under the
auspices of the group  GNFM-INdAM and is partially supported by MIUR via the PRIN 2017 ``Hyperbolic Systems of Conservation Laws and Fluid Dynamics: Analysis and Applications''.
 The research activity of F. Crispo  is also supported by GNFM-INdAM via Progetto Giovani 2017.}

\end{document}